# THE POLARISED PARTITION RELATION
# FOR ORDER TYPES

LUKAS DANIEL KLAUSNER AND THILO WEINERT

ABSTRACT. We analyse partitions of products with two ordered factors in two classes where both factors are countable or well-ordered and at least one of them is countable. This relates the partition properties of these products to cardinal characteristics of the continuum. We build on work by Erdős, Garti, Jones, Orr, Rado, Shelah and Szemerédi. In particular, we show that a theorem of Jones extends from the natural numbers to the rational ones but consistently extends only to three further equimorphism classes of countable orderings. This is made possible by applying a thirteen-year old theorem of Orr about embedding a given order into a sum of finite orders indexed over the given order.

## 1. INTRODUCTION

The partition calculus was introduced over six decades ago by Erdős and Rado in their seminal paper [ER56]. They introduced the ordinary partition relation which concern partitions of finite subsets of a set of a given size and the polarised partition relation which concerns partitions of finite subsets of products of sets of a given size. The notion of "size" here was mostly taken to refer to the cardinality of a set, but can easily be interpreted to refer to other notions of size, for example the order type of an ordered set. Assuming the axiom of choice, every set can be well-ordered and there naturally is the smallest ordinal which can be the order type of such a well-order. In this context, it is convenient to refer to this ordinal as the cardinality of the set in question and the analysis of partition relations for order types thereby naturally includes the one of partition relations for cardinality.
Partition relations for order types have been investigated in many papers. Chapter 7 of the book [Wil77] was devoted to them, and recently there has been renewed interest, cf. [Sch10, Sch12, Wei14, Hil16, CH17, Mer19, LHW19, IRW]. In these papers, however, the focus has been on the ordinary partition relation. To the authors' knowledge, relations for partitions of products of order types other than initial ordinals have only been studied in [EHM70, EHM71].
In this paper we aim to study polarised partition relations for order types. We limit the study to partitions of products with two factors such that both are countable or well-ordered and at least one of them is countable. In section 2, section 3 and section 4, we

2010 *Mathematics Subject Classification.* Primary 03E02; Secondary 03E17, 05D10, 06A05.
The first author was supported by the Austrian Science Fund (FWF) project P29575 "Forcing Methods: Creatures, Products and Iterations".
The second author was supported by the Austrian Science Fund (FWF) projects I3081 "Filters, Ultrafilters and Connections with Forcing" and Y1012 "Infinitary Combinatorics and Definability".
We are grateful to Raphaël Carroy, Garrett Ervin and Martin Goldstern for inspiring discussions.

briefly review the subjects of cardinal characteristics of the continuum, order types and polarised partition relations, respectively.

In section 5, we deal with the case of two countable sources and derive some results building on Ramsey's theorem and other results about the ordinary partition relation. Theorem 5.12 in particular establishes a positive partition relation for any two indecomposable ordinals smaller than $\omega^\omega$ as factors. In section 6, we deal with the case of one countable and one well-ordered source. Theorem 6.7 improves on a theorem of Jones by weakening its hypothesis and strengthening its conclusion – the latter by replacing the order type of the natural numbers by the one of the rational numbers. Theorem 6.17 limits the scope of further strengthenings in this manner as it establishes a negative partition relation involving the unbounding number $\mathfrak{b}$, a cardinal characteristic of the continuum. Such cardinal characteristics have been extensively studied over the past decades and many independence results regarding them have been attained. Via our results, some of these independence results directly lead to the independence of certain partition relations over ZFC. Connections between cardinal characteristics and partition relations have been explored before in [GS12, GS14, RT18, LHW19, CGW20]. It is worth noting that one of the difficulties of working with non-well-ordered types in the proof of Theorem 6.17 was resolved quite elegantly by employing a theorem of Orr. We conclude the paper by asking a few questions, some of which we expect to be decided within the framework of ZFC.

## 2. Cardinal Characteristics

We briefly recall the definitions of some cardinal characteristics relevant to our discussion. For further background, cf. [vD84, BJ95, Bla10, Hal17]. For functions $g, h \colon \omega \longrightarrow \omega$, we say that $g$ *eventually dominates* $h$ if the set of all natural numbers $n$ such that $g(n) \leqslant h(n)$ is finite. For sets $x$ and $y$ of natural numbers we say that $x$ *splits* $y$ if both $y \cap x$ and $y \setminus x$ are infinite. We say that $x$ *almost contains* $y$ if $y \setminus x$ is finite.

**Definition 2.1** ([vD84])**.** An *unbounded family* is a family $F$ of functions $g \colon \omega \longrightarrow \omega$ such that no single function $h \colon \omega \longrightarrow \omega$ eventually dominates all members of $F$. The *unbounding number* (sometimes called the bounding number) $\mathfrak{b}$ is the smallest cardinality of an unbounded family.

**Definition 2.2** ([vD84])**.** A *splitting family* is a family $F$ of sets of natural numbers such that for every infinite set $x$ of natural numbers, there is a member of $F$ splitting $x$. The *splitting number* $\mathfrak{s}$ is the smallest cardinality of a splitting family.

**Definition 2.3.** A *countably splitting family* is a family $F$ of sets of natural numbers such that for every countable collection $X$ of infinite sets of natural numbers, there is a member of $F$ splitting every element of $X$. The *countably splitting number* $\mathfrak{s}_{\aleph_0}$ is the smallest cardinality of a countably splitting family.

**Definition 2.4** ([vD84])**.** A *tower* is a sequence $\langle x_\xi \mid \xi < \alpha \rangle$ of infinite sets of natural numbers such that for $\gamma < \beta$, the set $x_\gamma$ almost contains $x_\beta$. A tower is *extendible* if there is an infinite set almost contained in every member of it. The *tower number* $\mathfrak{t}$ is the smallest ordinal $\alpha$ such that not all towers of length $\alpha$ are extendible.

Also recall that $\text{cov}(\mathcal{M})$ denotes the minimal number of meagre sets of reals necessary to cover the reals.

**Observation 2.5.** $\mathfrak{s} \leqslant \mathfrak{s}_{\aleph_0}$.

**Theorem 2.6.** $\mathfrak{s}_{\aleph_0} \leqslant \max(\mathfrak{b}, \mathfrak{s})$.

**Theorem 2.7.** $\min(\text{cov}(\mathcal{M}), \mathfrak{s}_{\aleph_0}) \leqslant \mathfrak{s}$.

Theorem 2.6 is due to Kamburelis and Węglorz, cf. [KW96, Proposition 2.1]. Theorem 2.7 is due to Kamburelis, cf. [KW96, Proposition 2.3]. It is still an unanswered question whether $\mathfrak{s} < \mathfrak{s}_{\aleph_0}$ is consistent. The notion of a countably splitting family was introduced by Malyhin in [Mal89]; also cf. [Ste93, Question 5.3]. We also note that the theorems we quote which mention the tower number were originally proved for the pseudointersection number $\mathfrak{p}$. By the recent seminal result by Malliaris and Shelah that $\mathfrak{p} = \mathfrak{t}$, the theorems follow in the form stated here. For the proof that the pseudointersection number equals the tower number, cf. [MS13].

As the unbounding number and the splitting number seem most relevant to the results in this paper, we would like to point out that their behaviour on regular cardinals was shown to be uncorrelated by Fischer et al., cf. [FS08, BF11, FM17].

For an overview of the relations between these cardinal characteristics, see Figure 1, where a line indicates a ZFC-provable inequality.

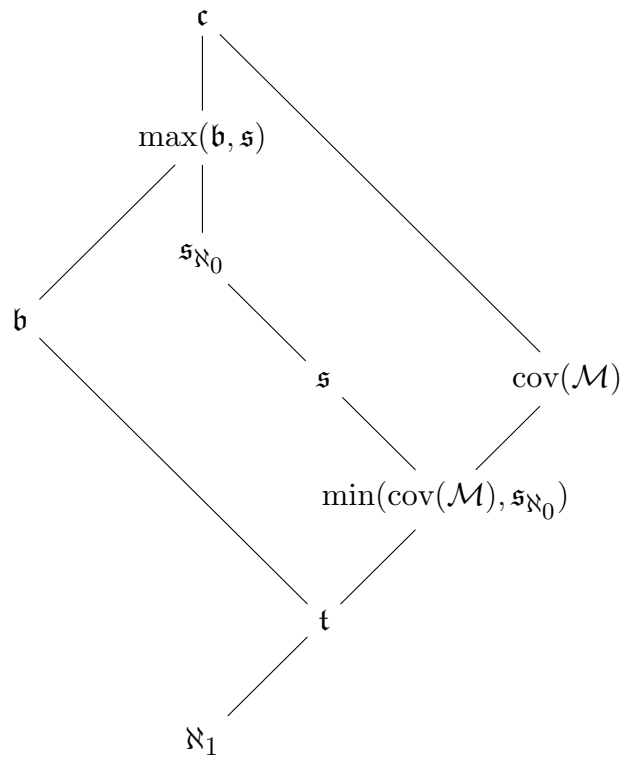

FIGURE 1. The inequalities between the aforementioned cardinal characteristics known to be ZFC-provable.

## 3. Order Types

By an *order type* we understand an equivalence class of ordered sets equivalent under order-preserving bijections. For an order type $\varphi$, we denote its reverse by $\varphi^*$. We denote the order type of the natural numbers by $\omega$ and the order type of the rational numbers by $\eta$. If $\varphi$ and $\psi$ are order types, $\varphi + \psi$ denotes the order type of the set $F \cup P$ under the ordering given by $(x \in F \wedge y \in P)$ or $(x <_F y)$ or $(x <_P y)$, where $F$ and $P$ are disjoint sets of order types $\varphi$ and $\psi$, respectively. Furthermore, $\varphi\psi$ denotes the order type of $P \times F$ under the (lexicographic) ordering $\langle p_0, f_0 \rangle < \langle p_1, f_1 \rangle$ given by $p_0 <_P p_1$ or $(p_0 = p_1 \wedge f_0 <_F f_1)$.

If $\varphi$ and $\psi$ are order types, we write $\varphi \leqslant \psi$ to indicate that an ordered set of type $\varphi$ may be embedded order-preservingly into an ordered set of type $\psi$. Contrary to the class of ordinals, which is well-ordered under embeddability, even the class of linear orders is only quasi-ordered by the embeddability relation. Embeddability is not a linear order between the linear orders, as shown by $\omega$ and $\omega^*$. It also fails to be antisymmetric, as witnessed by $\eta$ and $\eta + 1$. An ordering which fails to embed $\eta$ is called *scattered*. Orderings which are presentable as countable unions of scattered orderings are called $\sigma$-scattered. By a seminal result of Laver, cf. [Lav71], the class of $\sigma$-scattered orderings is well-quasi-ordered by the embeddability relation. Recall that an ordering is a *well-quasi-ordering* if both every strictly descending sequence and every set of pairwise incomparable elements is finite, cf. [Ros82, chapter 10] or [Fra00, chapter 4]. We say that two orders are *equimorphic* (sometimes also called *biembeddable*) if they may be embedded order-preservingly into each other. While equimorphic ordinals are identical, there are pairs of non-isomorphic yet equimorphic order types.

We call an order type $\varphi$ *additively decomposable* if there are types $\psi$ and $\tau$ such that $\varphi = \psi + \tau$ but neither $\varphi \leqslant \psi$ nor $\varphi \leqslant \tau$. We call it *unionwise decomposable* if there is an ordered set $\langle X, < \rangle$ of type $\varphi$ and a $Y \subseteq X$ such that neither $\varphi \leqslant \mathrm{otp}(\langle Y, < \rangle)$ nor $\varphi \leqslant \mathrm{otp}(\langle X \setminus Y, < \rangle)$. We call it *multiplicatively decomposable* if there are types $\psi$ and $\tau$ such that $\varphi = \psi\tau$ but neither $\varphi \leqslant \psi$ nor $\varphi \leqslant \tau$. We call it *typewise decomposable* if there is an ordered set $\langle X, <_X \rangle$ and for every $x \in X$ disjoint ordered sets $\langle Y_x, <_x \rangle$ such that the set $\langle \bigcup_{x \in X} Y_x, < \rangle$ has type $\varphi$ if $a < b$ is given by

$$\exists x \, (\exists y \colon a \in Y_x \wedge b \in Y_y \wedge x <_X y) \vee (a \in Y_x \wedge b \in Y_x \wedge a <_x b)$$

and furthermore neither $\varphi \leqslant \mathrm{otp}(\langle X, <_X \rangle)$ nor $\varphi \leqslant \mathrm{otp}(\langle Y_x, <_x \rangle)$ for any $x \in X$.

An order type is called (additively, unionwise, multiplicatively, typewise) *indecomposable* if it fails to be (additively, unionwise, multiplicatively, typewise) decomposable. (Unionwise indecomposability is called strong indecomposability in [EM72].) For ordinals, additive decomposability is equivalent to unionwise decomposability. An ordinal $\alpha$ is additively indecomposable if and only if there is an ordinal $\beta$ such that $\alpha = \omega^\beta$.

We call $\varphi$ *multiplicatively surpassable* if there are types $\psi$ and $\tau$ such that $\varphi \leqslant \psi\tau$ but neither $\varphi \leqslant \psi$ nor $\varphi \leqslant \tau$. We call $\varphi$ *multiplicatively unsurpassable* if it fails to be multiplicatively surpassable.

We call $\varphi$ *multiplicatively transcendable* if there are types $\psi$ and $\tau$ such that $\varphi \leqslant \psi\tau$ and $\psi \leqslant \varphi$ and $\tau \leqslant \varphi$ but neither $\varphi \leqslant \psi$ nor $\varphi \leqslant \tau$. We call $\varphi$ *multiplicatively untranscendable* if it fails to be multiplicatively transcendable.

Clearly, additive decomposability implies unionwise decomposability and multiplicative decomposability implies typewise decomposability.

Let us recall a famous Theorem of Hausdorff, cf. [Hau08, Satz XII].

**Theorem 3.1.** *The class of scattered order types is the smallest non-empty class containing all reversals and well-ordered sums.*

**Corollary 3.2.** *Up to equimorphism, the only countable typewise indecomposable order types are 0, 1, 2, $\omega$, $\omega^*$, and $\eta$.*

**Observation 3.3.**
- *Unionwise decomposability fails to imply additive decomposability as shown by the order type $\omega(\omega^* + \omega)$.*
- *Typewise decomposability fails to imply multiplicative decomposability as shown by the order type $\omega^\omega$.*
- *Multiplicative transcendability fails to imply multiplicative decomposability as shown by the order type $\omega + \omega^*$.*

## 4. Polarised Partition Relations

We will analyse the *polarised partition relation* for order types:

$$\begin{pmatrix} \alpha \\ \beta \end{pmatrix} \longrightarrow \begin{pmatrix} \gamma & \varepsilon \\ \delta & \zeta \end{pmatrix}.$$

This relation states that for every colouring $c\colon A \times B \longrightarrow 2$ of the product of a set $A$ of size $\alpha$ and a set $B$ of size $\beta$, either there is a $C \subseteq A$ of size $\gamma$ and a $D \subseteq B$ of size $\delta$ such that $c[C \times D] \subseteq \{0\}$ or there is an $E \subseteq A$ of size $\varepsilon$ and a $Z \subseteq D$ of size $\zeta$ such that $c[E \times Z] = \{1\}$.

The following relation states that for every ordinal $\kappa$ and every colouring $c\colon A \times B \longrightarrow \kappa$ of a set $A$ of size $\alpha$ and a set $B$ of size $\beta$, there is a $\vartheta < \kappa$, a $C \subseteq A$ of size $\gamma$ and a $D \subseteq B$ of size $\delta$ such that $c[C \times D] = \{\vartheta\}$:

$$\begin{pmatrix} \alpha \\ \beta \end{pmatrix} \longrightarrow \begin{pmatrix} \gamma \\ \delta \end{pmatrix}_\kappa.$$

The term "size" is intentionally vague in these statements and will be taken to mean cardinality or order type depending on the context. Note that

$$\begin{pmatrix} \alpha \\ \beta \end{pmatrix} \longrightarrow \begin{pmatrix} \gamma & \gamma \\ \delta & \delta \end{pmatrix}$$

is equivalent to

$$\begin{pmatrix} \alpha \\ \beta \end{pmatrix} \longrightarrow \begin{pmatrix} \gamma \\ \delta \end{pmatrix}_2.$$

It is common to call the parameters to the left of the arrow *sources* and the ones to the right of the arrow *targets* of the relation. The failure of the relations above is expressed by crossed arrows, i. e.

$$\begin{pmatrix}\alpha\\ \beta\end{pmatrix} \not\longrightarrow \begin{pmatrix}\gamma & \varepsilon\\ \delta & \zeta\end{pmatrix}$$

and

$$\begin{pmatrix}\alpha\\ \beta\end{pmatrix} \not\longrightarrow \begin{pmatrix}\gamma\\ \delta\end{pmatrix}_\varepsilon,$$

respectively.

We will also use the *ordinary partition relation*: $\alpha \longrightarrow (\beta_0, \ldots, \beta_k)^n$ states that for every colouring $c\colon [A]^n \longrightarrow k+1$ of the unordered $n$-tuples with $k+1$ colours, there is an $i \leqslant k$ and a $C \subseteq A$ of size $\beta_i$ such that $c\bigl[[C]^n\bigr] \subseteq \{i\}$. Moreover, $\alpha \longrightarrow (\beta)^n_m$ stands for $\alpha \longrightarrow (\underbrace{\beta, \ldots, \beta}_{m \text{ times}})^n$.

## 5. Both Sources Countable

The following fact limits the realm of relations holding for countable sources; it is easily verified by enumerating both sources in order type $\omega$ and colouring a pair of one element each of both sources by comparing their respective indices in these enumerations.

**Fact 5.1.**
$$\begin{pmatrix}\eta\\ \eta\end{pmatrix} \not\longrightarrow \begin{pmatrix}1 & \aleph_0\\ \aleph_0 & 1\end{pmatrix}.$$

The following fact can be verified by use of the pigeonhole principle and the notion of unionwise indecomposablity.

**Fact 5.2.** *For all natural numbers $m$, $n$ and all unionwise indecomposable types $\varphi$,*
$$\begin{pmatrix}\varphi\\ mn+1\end{pmatrix} \longrightarrow \begin{pmatrix}\varphi\\ n+1\end{pmatrix}_m.$$

Considering these facts, it is natural to ask for which countable linear order types $\varphi$, $\psi$ and natural numbers $n$ the relation

$$\begin{pmatrix}\varphi\\ \psi\end{pmatrix} \longrightarrow \begin{pmatrix}\varphi & n\\ \psi & n\end{pmatrix}$$

holds.

We can partially answer this using Ramsey's theorem, a technique which was used in [HS69] to analyse ordinary partition relations involving finite powers of $\omega$.

**Lemma 5.3.** *If $k$, $m$ and $n$ are natural numbers, then*
$$\begin{pmatrix}\omega^k\\ \omega^m\end{pmatrix} \longrightarrow \begin{pmatrix}\omega^k & n\\ \omega^m & n\end{pmatrix}.$$

*Proof.* Let $k$, $m$ and $n$ be natural numbers. Let
$$\left\langle \langle n_{i,j} \mid j < k \rangle \,\Big|\, i < \binom{k+m}{k} \right\rangle$$
be an enumeration of all descending enumerations of $[k+m]^k$, and for all $i < \binom{k+m}{k}$, let $\langle n_{i,j} \mid j \in k+m \setminus k \rangle$ be the descending enumeration of $k+m \setminus \{n_{i,j} \mid j < k\}$. Now let $c \colon \omega^k \times \omega^m \longrightarrow 2$. We may define a colouring $c' \colon [\omega]^{k+m} \longrightarrow 2^{\binom{k+m}{k}}$ by

$$(*_1) \quad \{\ell_i \mid i < k+m\}_< \longmapsto \sum_{i < \binom{k+m}{k}} 2^i \, c\Big(\omega^{k-1}\ell_{n_{i,k-1}} + \cdots + \ell_{n_{i,0}}, \omega^{m-1}\ell_{n_{i,k+m-1}} + \cdots + \ell_{n_{i,k}}\Big).$$

Now by Ramsey's theorem, there is an infinite set $H$ of natural numbers homogeneous for $c'$. We distinguish two cases.

First, assume that $H$ is homogeneous for $c'$ in colour 0. Consider a partition $\{H_0, H_1\}$ of $H$ into two infinite sets. Furthermore consider the sets
$$X := \{\omega^{k-1}n'_{k-1} + \cdots + n'_0 \mid \langle n'_i \mid i < k \rangle \text{ is descending and } \{n'_i \mid i < k\} \subseteq H_0\}$$
and
$$Y := \{\omega^{m-1}\tilde{n}_{m-1} + \cdots + \tilde{n}_0 \mid \langle \tilde{n}_i \mid i < m \rangle \text{ is descending and } \{\tilde{n}_i \mid i < m\} \subseteq H_1\}.$$

Note that $X$ has order type $\omega^k$ while $Y$ has order type $\omega^m$. We will prove that $X \times Y$ is homogeneous for $c$ in colour 0. To this end, let $\omega^{k-1}n'_{k-1} + \cdots + n'_0 = x \in X$ and $\omega^{m-1}\tilde{n}_{m-1} + \cdots + \tilde{n}_0 = y \in Y$. Consider the ascending enumeration $\langle \ell_i \mid i < k+m \rangle$ of $x \cup y$. There is an $i < \binom{k+m}{k}$ such that $n'_j = \ell_{n_{i,j}}$ for $j < k$ and $\tilde{n}_j = \ell_{n_{i,k+j}}$ for $j < m$. As $H$ is homogeneous for $c'$ in colour 0, we have $c'(x \cup y) = 0$ and as $c'$ assigns a set the value 0 if and only if all summands in Eq. $(*_1)$ have value 0, we also have $c(x,y) = 0$. This concludes the argument for $X \times Y$ being homogeneous for $c$ in colour 0.

Now assume that $H$ is homogeneous for $c'$ in a positive colour. This means that there is an $i < \binom{k+m}{k}$ such that $c\big(\omega^{k-1}\ell_{n_{i,k-1}} + \cdots + \ell_{n_{i,0}}, \omega^{m-1}\ell_{n_{i,k+m-1}} + \cdots + \ell_{n_{i,k}}\big) = 1$ for all $\{\ell_j \mid j < k+m\}_< \in [\omega]^{k+m}$. Now we are going to inductively construct two disjoint sets $X, Y \subseteq H$ of cardinality $n$ such that $c(x,y) = 1$ for all $x \in X$ and $y \in Y$. Let $\langle h_j \mid j < (k+m) \cdot n \rangle$ be an ascending enumeration of elements of $H$. We define $x_j := \omega^{k-1}h_{n_{i,k-1} \cdot n + j} + \cdots + h_{n_{i,0} \cdot n + j}$ and $y_j := \omega^{m-1}h_{n_{i,k+m-1} \cdot n + j} + \cdots + h_{n_{i,k} \cdot n + j}$ for $j < n$. We furthermore define $X := \{x_j \mid j < n\}$ and $Y := \{y_j \mid j < n\}$. It is easy to see that $c(x,y) = 1$ for all $x \in X$ and $y \in Y$. □

**Lemma 5.4.** *For all order types $\rho, \tau, \varphi$ and $\psi$, $\rho \longrightarrow (2\tau, \varphi + \psi, \psi + \varphi)^2$ implies*
$$\binom{\rho}{\rho} \longrightarrow \begin{pmatrix} \tau & \varphi \\ \tau & \psi \end{pmatrix}.$$

*Proof.* Let $\rho, \tau, \varphi$ and $\psi$ be order types as above and assume that $\rho \longrightarrow (2\tau, \varphi+\psi, \psi+\varphi)^2$. Furthermore, let $E \subseteq \rho \times \rho$. We define a colouring
$$c \colon [\rho]^2 \longrightarrow 3 \colon \{\nu, \xi\}_< \longmapsto [\langle \nu, \xi \rangle \in E] + 2[\langle \nu, \xi \rangle \notin E \ni \langle \xi, \nu \rangle]$$
As $\rho \longrightarrow (2\tau, \varphi+\psi, \psi+\varphi)^2$, we have to distinguish three cases.

For the first case, assume that there is an $X \in [\rho]^{2\tau}$ homogeneous for $c$ in colour 0. Then there is $A \in [X]^{\tau}$ such that $X \setminus A \in [X]^{\tau}$ and clearly $A \times (X \setminus A) \subseteq \rho \times \rho \setminus E$.

For the second case, assume that there is a $Y \in [\rho]^{\varphi+\psi}$ homogeneous for $c$ in colour 1. Let $B \in [Y]^{\varphi}$ be such that $B < Y \setminus B$ and $Y \setminus B \in [Y]^{\psi}$. Then $B \times (Y \setminus B) \subseteq E$.

For the final case, assume that there is a $Z \in [\rho]^{\psi+\varphi}$ homogeneous for $c$ in colour 2. Let $C \in [Z]^{\psi}$ be such that $C < Z \setminus C$ and $Z \setminus C \in [Z]^{\varphi}$. Then $(Z \setminus C) \times C \subseteq E$. □

**Theorem 5.5** ([ER56, Theorem 6]). $\eta \longrightarrow (\eta, \aleph_0)^2$.

**Theorem 5.6** ([Lar73]). *For all natural numbers $n$, $\omega^\omega \longrightarrow (\omega^\omega, n)^2$.*

**Proposition 5.7.** *For all natural numbers $k$,*
$$\binom{\eta}{\eta} \longrightarrow \binom{\eta \ \ k}{\eta \ \ k}.$$

*Proof.* Let $k$ be any natural number. By Theorem 5.5, we have that $\eta \longrightarrow (\eta, \aleph_0)^2$, so in particular we have $\eta \longrightarrow (\eta, \binom{4k-2}{2k-1})^2$. As $\binom{4k-2}{2k-1} \longrightarrow (2k)^2_2$, cf. [ES35], this implies $\eta \longrightarrow (\eta, 2k, 2k)^2$ and (as $2\eta$ may be embedded into $\eta$) also $\eta \longrightarrow (2\eta, 2k, 2k)^2$. Then by Lemma 5.4, the statement of the proposition follows. □

**Proposition 5.8.** *For all natural numbers $k$,*
$$\binom{\omega^\omega}{\omega^\omega} \longrightarrow \binom{\omega^\omega \ \ k}{\omega^\omega \ \ k}.$$

*Proof.* Let $k$ be any natural number. By Theorem 5.6, we have that $\omega^\omega \longrightarrow (\omega^\omega, \binom{4k-2}{2k-1})^2$. Via $\binom{4k-2}{2k-1} \longrightarrow (2k)^2_2$, cf. [ES35], this implies $\omega^\omega \longrightarrow (\omega^\omega, 2k, 2k)^2$. Note that $2\omega^\omega = \omega^\omega$, so Lemma 5.4 implies the statement of the proposition. □

At this point we would like to recall the notion of *pinning*, cf. [GL75].

**Definition 5.9.** An order type $\varphi$ can be *pinned* to an order type $\psi$ (written as $\varphi \to \psi$) if for every ordered set $F$ of type $\varphi$ and $P$ of type $\psi$ there is a function (a so-called *pinning map*) $f \colon F \longrightarrow P$ such that $f[X] \in [P]^{\psi}$ for every $X \in [F]^{\varphi}$.

Clearly, every order type can be pinned to the initial ordinal of its cardinality. If $\xi \to \rho$, then replacing all occurences of $\xi$ by $\rho$ in a valid partition relation yields a partition relation which itself again is valid. In particular:

$$(*_2) \qquad \binom{\sigma}{\xi} \longrightarrow \binom{\tau \ \ \varphi}{\xi \ \ \psi} \quad \text{implies} \quad \binom{\sigma}{\rho} \longrightarrow \binom{\tau \ \ \varphi}{\rho \ \ \psi}$$

and

$$(*_3) \qquad \binom{\sigma}{\xi} \longrightarrow \binom{\tau \ \ \varphi}{\xi \ \ \xi} \quad \text{implies} \quad \binom{\sigma}{\rho} \longrightarrow \binom{\tau \ \ \varphi}{\rho \ \ \rho}$$

for all order types $\xi, \rho, \sigma, \tau, \varphi, \psi$. Therefore, Eq. ($*_2$) implies that Proposition 5.7 and Proposition 5.8 have the following corollary.

**Corollary 5.10.** *For all natural numbers $k$,*
$$\binom{\eta}{\omega} \longrightarrow \binom{\eta \quad k}{\omega \quad k} \text{ and } \binom{\omega^\omega}{\omega} \longrightarrow \binom{\omega^\omega \quad k}{\omega \quad k}.$$

**Lemma 5.11.** *For all natural numbers $k$ and $m$ and all order types $\varphi$ and $\psi$ and collections of order types $\langle \sigma_i \mid i < k \rangle$ and $\langle \tau_j \mid j < m \rangle$, if*
$$\binom{\sigma_i}{\tau_j} \longrightarrow \binom{\sigma_i \quad \varphi}{\tau_j \quad \psi}$$
*for all $i < k$ and all $j < m$, then*
$$\binom{\sum_{i<k} \sigma_i}{\sum_{j<m} \tau_j} \longrightarrow \binom{\sum_{i<k} \sigma_i \quad \varphi}{\sum_{j<m} \tau_j \quad \psi}.$$

*Proof.* Assume towards a contradiction that
$$\binom{\sigma_i}{\tau_j} \longrightarrow \binom{\sigma_i \quad \varphi}{\tau_j \quad \psi}$$
for all $i < k$ and all $j < m$, but
$$\binom{\sum_{i<k} \sigma_i}{\sum_{j<m} \tau_j} \not\longrightarrow \binom{\sum_{i<k} \sigma_i \quad \varphi}{\sum_{j<m} \tau_j \quad \psi}.$$

Let $\langle \langle i_\ell, j_\ell \rangle \mid \ell < km \rangle$ be an enumeration of $k \times m$. Moreover, for $i < k$, let $X_i^{(0)}$ be a set ordered by $<_i$ of type $\sigma_i$ and for $j < m$, let $Y_j^{(0)}$ be a set ordered by $<_{k+j}$ of type $\tau_j$. Then $X := \bigcup_{i<k} X_i^{(0)}$ ordered by
$$\bigcup_{i<k} (<_i \cup \bigcup_{j<i} (X_j^{(0)} \times X_i^{(0)}))$$
has type $\sum_{i<k} \sigma_i$ and $Y := \bigcup_{i<m} Y_i^{(0)}$ ordered by
$$\bigcup_{i<m} (<_{k+i} \cup \bigcup_{j<i} (Y_j^{(0)} \times Y_i^{(0)}))$$
has type $\sum_{j<m} \tau_j$. As
$$(*_4) \qquad \binom{\sum_{i<k} \sigma_i}{\sum_{j<m} \tau_j} \not\longrightarrow \binom{\sum_{i<k} \sigma_i \quad \varphi}{\sum_{j<m} \tau_j \quad \psi},$$
we know that there is an $E \subseteq X \times Y$ witnessing this; in particular this $E$ is such that for no $A \in [X]^\varphi$ there is a $B \in [Y]^\psi$ such that $A \times B \subseteq E$. As
$$\binom{\sigma_i}{\tau_j} \longrightarrow \binom{\sigma_i \quad \varphi}{\tau_j \quad \psi}$$
for all $i < k$ and all $j < m$, we know that for all $i < k$, all $j < m$, all $A \in [X]^{\sigma_i}$ and all $B \in [Y]^{\tau_j}$, there are $C \in [A]^{\sigma_i}$ and $D \in [B]^{\tau_j}$ such that
$$(*_5) \qquad C \times D \subseteq X \times Y \setminus E.$$

Now suppose that we are in step $\ell$ of the induction and have by now defined $X_i^{(\ell)} \in [X_i^{(0)}]^{\sigma_i}$ for $i < k$ and $Y_j^{(\ell)} \in [Y_j^{(0)}]^{\tau_j}$ for $j < m$. Now Eq. ($*_5$) implies that there are $X_{i_\ell}^{(\ell+1)} \in [X_{i_\ell}^{(\ell)}]^{\sigma_{i_\ell}}$ and $Y_{j_\ell}^{(\ell+1)} \in [Y_{j_\ell}^{(\ell)}]^{\tau_{j_\ell}}$ such that $X_{i_\ell}^{(\ell+1)} \times Y_{j_\ell}^{(\ell+1)} \subseteq X \times Y \setminus E$. Let $X_i^{(\ell+1)} := X_i^{(\ell)}$ for all $i \in k \setminus \{i_\ell\}$ and $Y_j^{(\ell+1)} := Y_j^{(\ell)}$ for all $j \in m \setminus \{j_\ell\}$. After the induction, we have that $\bigcup_{i<k} X_i^{(km)} \subseteq X$ has order type $\sum_{i<k} \sigma_i$ and $\bigcup_{j<m} Y_j^{(km)}$ has order type $\sum_{j<m} \tau_j$. Moreover, $(\bigcup_{i<k} X_i^{(km)} \subseteq X) \times (\bigcup_{j<m} Y_j^{(km)}) \subseteq X \times Y \setminus E$. This contradicts $E$ being a witness to Eq. ($*_4$). □

**Theorem 5.12.** *For all ordinals $\alpha, \beta < \omega^\omega$ and all natural numbers $n$,*

$$\begin{pmatrix} \omega\alpha \\ \omega\beta \end{pmatrix} \longrightarrow \begin{pmatrix} \omega\alpha & n \\ \omega\beta & n \end{pmatrix}.$$

*Proof.* Note that if $\nu < \omega^\omega$, then $\omega\nu$ is a finite sum of infinite indecomposable ordinals smaller than $\omega^\omega$. So as $\{\alpha, \beta\} \subseteq \omega^\omega$, both $\omega\alpha$ and $\omega\beta$ are finite sums of infinite indecomposable ordinals smaller than $\omega^\omega$. Using Lemma 5.3 together with Lemma 5.11 proves the theorem. □

## 6. One Source Countable

Having gained some understanding of the case with two countable sources, we now consider the case of one countable and one well-ordered source. The guiding question is the following:

**Question A.** *For which countable linear order types $\varphi$ and which ordinals $\alpha, \beta$ do we have*

$$\begin{pmatrix} \omega_1 \\ \varphi \end{pmatrix} \longrightarrow \begin{pmatrix} \alpha & \beta \\ \varphi & \varphi \end{pmatrix}?$$

Let us note that we cannot hope for a relation in the case that $\varphi$ is decomposable.

**Observation 6.1.** *If $\varphi$ is a unionwise decomposable order type and $\psi$ is any order type, then*

$$\begin{pmatrix} \psi \\ \varphi \end{pmatrix} \nrightarrow \begin{pmatrix} 1 & 1 \\ \varphi & \varphi \end{pmatrix}.$$

The following is a corollary of a claim proved by Garti and Shelah in [GS14, Claim 1.2]:

**Corollary 6.2.** *If $\kappa < \mathfrak{s}$ is an infinite cardinal, then*

$$\begin{pmatrix} \kappa \\ \omega \end{pmatrix} \longrightarrow \begin{pmatrix} \kappa \\ \omega \end{pmatrix}_2 \text{ if and only if } \omega < \operatorname{cf}(\kappa).$$

We would like to point out an observation made by Brendle and Raghavan:

**Definition 6.3** ([BR14, Definition 31])**.** Let $\kappa$ be a regular cardinal, and let $\bar{A} = \langle a_\alpha \mid \alpha < \kappa \rangle$. $\bar{A}$ is *tail-splitting* if for every $b \in [\omega]^\omega$, there is $\alpha < \kappa$ such that $a_\beta$ splits $b$ for all $\beta \geqslant \alpha$. $\bar{A}$ is *club-splitting* if for every $b \in [\omega]^\omega$, $C_b = \{\alpha < \kappa \mid a_\alpha \text{ splits } b\}$ contains a club.

**Observation 6.4** ([BR14, Observation 44]). *The following are equivalent (for infinite ordinals $\mu$):*

- $\binom{\mu}{\omega} \longrightarrow \binom{\mu}{\omega}_2$
- $\mathrm{cf}(\mu) \neq \omega$ *and there does not exist a tail-splitting sequence of length $\mu$.*

As every tail-splitting sequence of a length with uncountable cofinality is club-splitting and every club-splitting sequence of a length with uncountable cofinality is countably splitting, we get the following corollary:

**Corollary 6.5.** *If $\kappa < \mathfrak{s}_{\aleph_0}$ is an infinite cardinal of uncountable cofinality, then*

$$\binom{\kappa}{\omega} \longrightarrow \binom{\kappa}{\omega}_2.$$

In [Jon08], Jones proved the following theorem, thus generalising an unpublished result of Szemerédi who proved it for the special case where $\kappa = \mathfrak{c}$, $\alpha$ is a cardinal and Martin's Axiom holds. (Szemerédi's result in turn generalised a result of Erdős and Rado from [ER56, penultimate page] who proved it for $\kappa = \omega_1$ and $\alpha = \omega$.)

**Theorem 6.6.** *For any regular uncountable $\kappa \leqslant \mathfrak{c}$ and any ordinal $\alpha < \min\{\mathfrak{t}, \kappa\}$,*

$$\binom{\kappa}{\omega} \longrightarrow \binom{\kappa \ \ \alpha}{\omega \ \ \omega}.$$

We can improve Jones' theorem by weakening its hypothesis and replacing the order type of the natural numbers $\omega$ by the order type of the rational numbers $\eta$; using pinning, this can be viewed as a strengthening of Jones' result.

**Theorem 6.7.** *For all cardinals $\kappa \leqslant \mathfrak{c}$ of uncountable cofinality and any ordinal $\alpha < \min\{\mathfrak{t}, \mathrm{cf}(\kappa)\}$,*

$$\binom{\kappa}{\eta} \longrightarrow \binom{\kappa \ \ \alpha}{\eta \ \ \eta}.$$

By pinning (cf. Definition 5.9), more precisely because of Eq. ($*_3$), we can immediately attain the following corollary.

**Corollary 6.8.** *For all cardinals $\kappa \leqslant \mathfrak{c}$ of uncountable cofinality, any ordinal $\alpha < \min\{\mathfrak{t}, \mathrm{cf}(\kappa)\}$ and both $\varphi \in \{\omega, \omega^*\}$,*

$$\binom{\kappa}{\varphi} \longrightarrow \binom{\kappa \ \ \alpha}{\varphi \ \ \varphi}.$$

First, note that by work of Bell, cf. [Bel81], its generalisation to all uncountable cardinals below the continuum, cf. [Bla10], and the recent breakthrough of Malliaris and Shelah, cf. [MS13], we have the following theorem.

**Theorem 6.9.** *For every cardinal $\kappa$, the following are equivalent:*

- $\kappa < \mathfrak{t}$.

- *Every family $F$ of cardinality $\kappa$ of infinite sets of natural numbers for which the intersection of finitely many members of $F$ is always infinite has an infinite pseudointersection. (Equivalently, $\kappa < \mathfrak{p}$.)*
- *Martin's Axiom holds for every family of cardinality $\kappa$ of dense subsets of a $\sigma$-centred partially ordered set. (Equivalently: For every $\sigma$-centred partially ordered set $P$ and every family $C$ of cardinality $\kappa$ of dense subsets of $P$, there is a filter on $P$ meeting every member of $C$.)*

Similar to Jones' proof, we require two lemmata. First, however, we make the following definition.

**Definition 6.10.** A family $\{A_\xi \mid \xi \in X\} \subseteq \mathcal{P}(\mathbb{Q})$ is an $\omega$-*complete filter basis dense in $\mathbb{Q}$* (or $\omega$-*complete $\mathbb{Q}$-basis* for short) if for any $x \in [X]^{<\omega}$, $\bigcap_{\xi \in x} A_\xi$ is dense in $\mathbb{Q}$.

**Lemma 6.11.** *Suppose $\kappa$ is a cardinal of uncountable cofinality. For $A \subseteq \mathbb{Q}$, let $\overline{A} := \mathbb{Q}\setminus A$. Then for any $\{A_\xi \mid \xi < \kappa\} \subseteq \mathcal{P}(\mathbb{Q})$ either*

*(a) there is $X \in [\kappa]^{\mathrm{cf}(\kappa)}$ such that $\{A_\xi \mid \xi \in X\}$ is an $\omega$-complete $\mathbb{Q}$-basis, or*
*(b) there is $Y \in [\kappa]^\kappa$ such that $\bigcap_{\xi \in Y} \overline{A_\xi}$ contains a copy of $\mathbb{Q}$.*

*Proof.* For $x \in [\kappa]^{<\omega}$, let $A(x) := \bigcap_{\xi \in x} A_\xi$. Suppose that (a) is false; we will show that (b) must hold. Fix $X \subseteq \kappa$ maximal such that $\{A_\xi \mid \xi \in X\}$ is an $\omega$-complete $\mathbb{Q}$-basis. Clearly, $|X| < \mathrm{cf}(\kappa)$.

For any $\delta \notin X$, there must be a set of witnesses $x_\delta \in [X]^{<\omega}$ and $r_\delta, s_\delta \in \mathbb{Q}$ showing the maximality of $X$, i.e. witnessing that $A_\delta \cap A(x_\delta)$ avoids the interval $(r_\delta, s_\delta)$ and hence is not dense in $\mathbb{Q}$.

Since $\mathrm{cf}(\kappa) > \aleph_0 = |\mathbb{Q}|$, by applying the pigeonhole principle, we can find some $Y$ of size $\kappa$ and some $x^* \in [X]^{<\omega}$ and $r^*, s^* \in \mathbb{Q}$ such that for any $\delta \in Y$, we have $x_\delta = x^*$ and $(r_\delta, s_\delta) = (r^*, s^*)$. Hence $\bigcup_{\delta \in Y} A_\delta$ avoids $A(x^*) \cap (r^*, s^*)$, which implies that $\bigcap_{\delta \in Y} \overline{A_\delta}$ contains $A(x^*) \cap (r^*, s^*)$ – which, by the definition of $X$, contains a copy of $\mathbb{Q}$. $\square$

While the proof of the previous lemma is a significant simplification of Jones' proof, the following lemma's proof is more or less the same as Jones'.

**Lemma 6.12.** *Suppose $\kappa$ is a cardinal of uncountable cofinality. Then for any $\omega$-complete $\mathbb{Q}$-basis $\{A_\xi \mid \xi < \kappa\} \subseteq \mathcal{P}(\mathbb{Q})$ and any ordinal $\alpha < \min\{\mathfrak{t}, \mathrm{cf}(\kappa)\}$, there is an $X \in [\kappa]^\alpha$ such that $\bigcap_{\xi \in X} A_\xi$ contains a copy of $\mathbb{Q}$.*

*Proof.* Let $\mu := |\omega + \alpha|$. Construct a sequence $\{M_\gamma \mid \gamma \leqslant \alpha\}$ of elementary submodels of $W := H_{\kappa^+}$ such that $|M_\gamma| = \mu$, $\{A_\xi \mid \xi < \kappa\} \in M_\gamma$, $\mu \cup \{\mu\} \subseteq M_\gamma$, and $M_\beta \in M_\gamma$ whenever $\beta < \gamma \leqslant \alpha$. Let $\vartheta := \sup(\kappa \cap M_\alpha)$ and note that $\vartheta < \kappa$ because $|M_\alpha| = \mu = |\omega + \alpha| < \mathrm{cf}(\kappa)$. Choose $\delta \in \kappa \setminus \vartheta$ arbitrarily.

For $x \in [\kappa]^{<\omega}$, let $A(x) := \bigcap_{\xi \in x} A_\xi$. For $a \in [\mathbb{Q}]^{<\omega}$, let $X(a) := \{\xi < \kappa \mid a \subseteq A_\xi\}$. Let $\mathbb{P}$ be the collection of all $p := \langle a_p, x_p \rangle$ with $a_p \in [A_\delta]^{<\omega}$, $x_p \in [\kappa \cap M_\alpha]^{<\omega}$ and $x_p \subseteq X(a_p)$ (or equivalently, $a_p \subseteq A(x_p)$). Write $q \leqslant p$ ("$q$ is stronger than $p$") if $a_p \subseteq a_q$ and $x_p \subseteq x_q$ (i.e. stronger conditions consist of larger sets), and note that $\mathbb{P}$ is $\sigma$-centred.

For $r, s \in \mathbb{Q}$, let $D_{r,s} := \{p \in \mathbb{P} \mid a_p \not\subseteq \mathbb{Q} \setminus (r, s)\}$. Each such $D_{r,s}$ is dense in $\mathbb{P}$, since given any $p \in \mathbb{P}$, $A(x_p)$ is dense in $\mathbb{Q}$, and hence we can pick some $t \in A(x_p) \cap (r, s)$ and define $q := \langle a_p \cup \{t\}, x_p \rangle$, which is stronger than $p$ and an element of $D_{r,s}$.

For $\xi < \alpha$, let $E_\xi := \{p \in \mathbb{P} \mid x_p \cap M_{\xi+1} \not\subseteq M_\xi\}$. Each such $E_\xi$ is dense in $\mathbb{P}$, since for any $p = \langle a_p, x_p \rangle \in \mathbb{P}$, the following consideration holds: In $W$, the formula
$$\forall \nu < \kappa^W = |\{A_\xi \mid \xi < \kappa\}| \; \exists \gamma > \nu \colon a_p \subseteq A_\gamma$$
is true (namely, one can pick $\gamma := \delta$); by elementarity, the formula must also hold in $M_{\xi+1}$ for any $\nu < \kappa \cap M_{\xi+1}$ – specifically for $\nu := \kappa \cap M_\xi$, and so we can find some $\gamma > \nu$ in $M_{\xi+1}$ such that $a_p \subseteq A_\gamma$. Letting $q := \langle a_p, x_p \cup \{\gamma\} \rangle$, we have found a condition $q$ which is stronger than $p$ and an element of $E_\xi$.

Since $|\omega + \alpha| < \mathfrak{t}$, there is a $\mathbb{P}$-generic filter $G$ meeting all $D_{r,s}$ and $E_\xi$ simultaneously. Letting $A := \bigcup_{p \in G} a_p$ and $X := \bigcup_{p \in G} x_p$, $A \subseteq \bigcap_{\xi \in X} A_\xi$ is dense in $\mathbb{Q}$ (and hence $\bigcap_{\xi \in X} A_\xi$ contains a copy of $\mathbb{Q}$) and $\mathrm{otp}(X) \geqslant \alpha$. $\square$

*Proof of Theorem 6.7.* Let $\mathbb{Q} \times \kappa =: K_0 \cup K_1$. For each $\xi < \kappa$ and $i \in \{0, 1\}$, let $K_i^\xi := \{r \in \mathbb{Q} \mid \langle r, \xi \rangle \in K_i\}$. By Lemma 6.11, there is either (a) an $X \in [\kappa]^{\mathrm{cf}(\kappa)}$ such that $\{K_1^\xi \mid \xi \in X\}$ is an $\omega$-complete $\mathbb{Q}$-basis or (b) a $Y \in [\kappa]^\kappa$ such that $\bigcap_{\xi \in Y} K_0^\xi$ contains a copy $Q_0$ of $\mathbb{Q}$. In case (b), simply let $A_0 := Y \in [\kappa]^\kappa$ and $B_0 := Q_0 \in [\mathbb{Q}]^\eta$; then $B_0 \times A_0 \subseteq K_0$. In case (a), by Lemma 6.12 there is an $X' \in [X]^\alpha$ such that $\bigcap_{\xi \in A_1} K_1^\xi$ contains a copy $Q_1$ of $\mathbb{Q}$. Let $A_1 := X' \in [\kappa]^\alpha$ and $B_1 := Q_1 \in [\mathbb{Q}]^\eta$; then $B_1 \times A_1 \subseteq K_1$. $\square$

Now we are going to show that consistently, other order types can take the place of $\omega$ and $\eta$ in the preceding results.

**Proposition 6.13.** *If $\kappa < \mathfrak{b}$ is a cardinal of uncountable cofinality while $n$ is a natural number and $\alpha \leqslant \kappa$, then*
$$\binom{\kappa}{\omega^n} \longrightarrow \begin{pmatrix} \kappa & \alpha \\ \omega^n & \omega^n \end{pmatrix} \text{ if and only if } \binom{\kappa}{\omega} \longrightarrow \begin{pmatrix} \kappa & \alpha \\ \omega & \omega \end{pmatrix}.$$

*Proof.* We will prove the theorem by induction. It obviously holds in case $n < 2$. So assume that $\kappa < \mathfrak{b}$ is a cardinal of uncountable cofinality, that $n$ is a natural number, that $\alpha \leqslant \kappa$ and that the equivalence has been proved up to $n$. As for any positive ordinal $n$, the ordinal $\omega^n$ may be pinned (cf. Definition 5.9) to $\omega$, the first statement implies the second. In order to show the converse, assume that the second statement holds.

Towards a contradiction, let $c \colon \kappa \times \omega^{n+1} \longrightarrow 2$ be such that $1 \in c[A \times B]$ for every $A \in [\kappa]^\kappa$ and $B \in [\omega^{n+1}]^{\omega^{n+1}}$ while $0 \in c[A \times B]$ for every $A \in [\kappa]^\alpha$ and $B \in [\omega^{n+1}]^{\omega^{n+1}}$.

Let $b \colon \omega \longrightarrow \omega^n$ be a bijection. For every countable ordinal $\xi$, we are going to define an infinite set $X_\xi$ of natural numbers and a sequence $\langle i_{\xi,\nu} \mid \nu < \omega^n \rangle$ such that for all ordinals $\nu < \omega^n$, the set $\{n \in X_\xi \mid c(\xi, \omega\nu + n) \neq i_{\xi,\nu}\}$ is finite. Let $X'_{\xi,0} \in [\omega]^\omega$ and $i_{\xi,b(0)}$ be such that $c(\xi, \omega b(0) + \ell) = i_{\xi,b(0)}$ for all $\xi < \kappa$ and all $\ell \in X'_{\xi,0}$. Suppose that there is a natural number $k$ such that we have defined $X'_{\xi,\nu}$ and $i_{\xi,\nu}$ for all $\nu < \omega^n$ such that $b^{-1}(\nu) \leqslant k$. We may choose $X'_{\xi,b(k+1)} \in [X'_{\xi,b(k)}]^\omega$ and $i_{\xi,b(k+1)}$ such that $c(\xi, \omega b(k+1) + \ell) = i_{\xi,b(k+1)}$ for all $\xi < \kappa$ and all $\ell \in X'_{\xi,b(k+1)}$. For every $\xi < \kappa$, let $X_\xi$ be a pseudointersection of the $X'_{\xi,\nu}$ for $\nu < \omega^n$.

Now we define a family $\{f_\xi \mid \xi < \kappa\}$ of functions

$$f_\xi \colon \omega \longrightarrow \omega \colon k \longmapsto \max \{n \in X_\xi \mid c(\xi, \omega b(k) + n) \neq i_{\xi,k}\}.$$

Now, as $\mathfrak{b} > \kappa$, there is a $g \colon \omega \longrightarrow \omega$ which eventually properly dominates every member of $\{f_\xi \mid \xi < \kappa\}$. For $\xi < \kappa$, let $n_\xi$ be a natural number such that $g$ dominates $f_\xi$ from $n_\xi$ on. As $\kappa$ has uncountable cofinality, there is a natural number $n$ and an $A \in [\kappa]^\kappa$ such that $n_\xi = n$ for all $\xi \in A$. We define a colouring

$$\psi \colon \kappa \times \omega \longrightarrow 2 \colon \langle \xi, k \rangle \longmapsto i_{\xi,k},$$

and as we assumed

$$\binom{\kappa}{\omega} \longrightarrow \binom{\kappa \ \ \alpha}{\omega \ \ \omega},$$

we may consider two cases.

First, assume that there is $A' \in [A]^\kappa$ and $Z \in [\omega]^\omega$ such that $\psi[A' \times Z] = \{0\}$. We consider the set $B := \{\omega b(k) + m \mid k \in Z \setminus n \ \wedge \ m \in \omega \setminus g(k)\}$. Clearly, $B$ has order type $\omega^{n+1}$. There has to be a $\xi \in A'$ and an $\nu \in B$ such that $c(\xi, \nu) = 1$. There is a $k \in Z \setminus n$ and an $m \in \omega \setminus g(k)$ such that $\nu = \omega b(k) + m$. As $c(\xi, \omega b(k) + m) = 1 \neq 0 = i_{\xi,k} = \psi(\xi, k)$, we have $m \leqslant f_\xi(k)$. But $g$ properly dominates $f_\xi$ from $n$ on, so $m \leqslant f_\xi(k) < g(k) \leqslant m$, a contradiction.

The second case works completely analogously. We assume that there is an $A' \in [A]^\alpha$ and a $Z \in [\omega]^\omega$ such that $\psi[A' \times Z] = \{1\}$. We consider the set $B := \{\omega b(k) + m \mid k \in Z \setminus n \wedge m \in \omega \setminus g(k)\}$. Clearly, $B$ has order type $\omega^{n+1}$. There has to be a $\xi \in A'$ and a $\nu \in B$ such that $c(\xi, \nu) = 0$. There is a $k \in Z \setminus n$ and an $m \in \omega \setminus g(k)$ such that $\nu = \omega b(k) + m$. As $c(\xi, \omega b(k) + m) = 0 \neq 1 = i_{\xi,k} = \psi(\xi, k)$, we have $m \leqslant f_\xi(k)$. But $g$ properly dominates $f_\xi$ from $n$ on, so $m \leqslant f_\xi(k) < g(k) \leqslant m$, a contradiction. $\square$

By Theorem 6.6 we get the following.

**Corollary 6.14.** *If $\kappa$ is a regular uncountable cardinal smaller than $\mathfrak{b}$ while $\beta \in \omega^\omega \setminus \omega$ is additively indecomposable and $\alpha < \min(\mathfrak{t}, \kappa)$, then*

$$\binom{\kappa}{\beta} \longrightarrow \binom{\kappa \ \ \alpha}{\beta \ \ \beta}.$$

Corollary 6.5 immediately yields another corollary:

**Corollary 6.15.** *If $\kappa < \min(\mathfrak{b}, \mathfrak{s}_{\aleph_0})$ is an infinite cardinal and $\beta \in \omega^\omega \setminus \omega$ is additively indecomposable, then*

$$\binom{\kappa}{\beta} \longrightarrow \binom{\kappa}{\beta}_2 \quad \text{if and only if } \omega < \operatorname{cf}(\kappa).$$

In order to prove the last theorem we reference a result of Orr.

**Proposition 6.16** ([Orr95, Proposition 2]). *Let $A$ be a countable linearly ordered set and for every $a \in A$ let $L_a$ be a finite linearly ordered set. Then there is an increasing map*

$$\sigma\colon A \longrightarrow L = \sum_{a \in A} L_a$$

*which maps onto all but finitely many points of $L$, and, in any event, onto at least one point in every $L_a$.*

Our final theorem shows that the assumption of $\kappa < \mathfrak{b}$ in Proposition 6.13 was indeed necessary.

**Theorem 6.17.** *If $\alpha$ is an ordinal of cofinality $\mathfrak{b}$ and $\varphi$ is a countable typewise decomposable order type, then*

$$\begin{pmatrix}\alpha\\\varphi\end{pmatrix} \not\rightarrow \begin{pmatrix}\alpha & 1\\\varphi & \varphi\end{pmatrix}.$$

*Proof.* Let $\alpha$ be an ordinal of cofinality $\mathfrak{b}$ and let this be witnessed by an ascending sequence $\langle \mu_\xi \mid \xi < \mathfrak{b}\rangle$ of ordinals cofinal in $\alpha$. Let $\langle f_\xi \mid \xi < \mathfrak{b}\rangle$ be an unbounded sequence of functions. We may assume without loss of generality that it is increasing modulo finite and every element of it is itself increasing. Let $\varphi$ be a countable typewise decomposable order type. Considering Observation 6.1, we may assume without loss of generality that $\varphi$ is unionwise indecomposable. As $\varphi$ is typewise decomposable, we have $3 \leqslant \varphi$, so as $\varphi$ is unionwise indecomposable, we have that $\varphi$ is infinite. Now for every $\xi < \alpha$, let $\nu_\xi$ be the smallest $\nu < \mathfrak{b}$ such that $\mu_\nu \geqslant \xi$. Furthermore let $\langle T, <_T\rangle$ be an ordered set and for every $t \in T$ let $\langle P_t, <_t\rangle$ be an ordered set such that $\sum_{t \in T} P_t$ has order type $\varphi$ and so witnesses the typewise decomposability of $\varphi$, i.e. neither $\varphi \leqslant \mathrm{otp}(T)$ nor $\varphi \leqslant P_t$ for any $t \in T$. Moreover, let $b\colon \omega \longrightarrow T$ and, for every natural number $n$, let $e_n\colon \omega \longrightarrow P_{b(n)}$ be enumerations. We define

$$(*_6) \qquad E := \{\langle \beta, \langle p, t\rangle\rangle \mid \beta < \alpha \wedge t \in T \wedge p \in P_t \wedge f_{\nu_\beta}(b^{-1}(t)) > e_{b^{-1}(t)}^{-1}(p)\}.$$

For the first half of the proof, let $X \in [\alpha]^\alpha$ and $Y \in [\sum_{t \in T} P_t]^\varphi$. We assume towards a contradiction that $X \times Y \subseteq \alpha \times \sum_{t \in T} P_t \setminus E$ and define

$$(*_7) \qquad a\colon \omega \longrightarrow \omega\colon n \longmapsto \min\{k \in \omega \setminus n \mid |Y \cap P_{b(k)}| = \aleph_0\}.$$

We first have to make sure that $a$ is well-defined.

**Claim 1.** *Eq. $(*_7)$ defines a function.*

*Proof of Claim.* Assume towards a contradiction that for all but finitely many $t \in T$ the set $Y \cap P_t$ is finite. Let $T' := \{t \in T \mid Y \cap P_t \text{ is finite}\}$. Then, by Proposition 6.16, there is a finite $F \subseteq \sum_{t \in T'}(Y \cap P_t)$ such that there is an increasing surjection $\sigma\colon T' \longrightarrow \sum_{t \in T'}(Y \cap P_t) \setminus F$. Using the Axiom of Choice, we may invert $\sigma$ to find that $\sum_{t \in T'}(Y \cap P_t) \setminus F \leqslant T'$. As $\varphi \not\leqslant \mathrm{otp}(T')$, we have $\varphi \not\leqslant \sum_{t \in T'}(Y \cap P_t) \setminus F$. As $\varphi$ is unionwise indecomposable and infinite and $F$ is finite, it follows that $\varphi \not\leqslant \sum_{t \in T'}(Y \cap P_t)$. Again, as $\varphi$ is unionwise indecomposable (and $\mathrm{otp}(Y) = \varphi$), we have $\varphi \leqslant \sum_{t \in T \setminus T'}(Y \cap P_t)$. Let $\{t_k \mid k < n\}$ be the ascending enumeration of $T \setminus T'$. Let $m < n$ be minimal such that $\varphi \leqslant \sum_{k \leqslant m}(Y \cap P_{t_k})$.

As $\varphi \not\leqslant P_{t_m}$, we have $\varphi \not\leqslant Y \cap P_{t_m}$, and as $\varphi$ is unionwise indecomposable, it follows that $\varphi \leqslant \sum_{k<m} P_{t_k}$, contradicting the minimality of $m$ and thus concluding the proof of the claim. □

We now define

$$(*_8) \quad \begin{aligned} g\colon \omega &\longrightarrow \omega\colon n \longmapsto \min\{k < \omega \mid Y \cap P_{b(n)} = \varnothing \vee e_n(k) \in Y\} \text{ and} \\ h\colon \omega &\longrightarrow \omega\colon n \longmapsto \max\{a(n), g(a(n))\}. \end{aligned}$$

As $\langle f_\xi \mid \xi < \mathfrak{b} \rangle$ is unbounded and increasing modulo finite and $X$ has order type $\alpha$, we can find a $\beta \in X$ such that $f_{\nu_\beta}$ is not eventually bounded by $h$. Let $k$ be a natural number such that $f_{\nu_\beta}(k) > h(k)$. By Eq. $(*_7)$, the set $Y \cap P_{a(k)}$ is infinite. By Eq. $(*_8)$, $e_{a(k)}\bigl(g(a(k))\bigr) \in Y$. By Eq. $(*_6)$ and since $X \times Y \subseteq \alpha \times \sum_{t \in T} P_t \setminus E$, we have $f_{\nu_\beta}(a(k)) \leqslant g(a(k))$. Ultimately, we thus get

$$g(a(k)) \leqslant h(k) < f_{\nu_\beta}(k) \leqslant f_{\nu_\beta}(a(k)) \leqslant g(a(k)),$$

a contradiction.

For the second half of the proof, let $\gamma < \alpha$ and $Y \in [\sum_{t \in T} P_t]^\varphi$ and assume towards a contradiction that $\{\gamma\} \times Y \subseteq E$. We distinguish two cases.

For the first case, we assume that there is a $t \in T$ such that $Q := \{p \in P_t \mid \langle t, p \rangle \in Y\}$ is infinite. Let $q \in Q$ be such that $e_{b^{-1}(t)}^{-1}(q) \in \omega \setminus f_\gamma(b^{-1}(t))$. Then by Eq. $(*_6)$ we have $\langle \gamma, \langle t, p \rangle \rangle \notin E$, a contradiction.

For the second case, we assume that for all $t \in T$ the set $\{p \in P_t \mid \langle t, p \rangle \in Y\}$ is finite. Then by Proposition 6.16 there is a finite $F \subseteq Y$ and an increasing map $\sigma\colon T \longrightarrow Y \setminus F$ which is onto. Using the Axiom of Choice we may invert $\sigma$ to see that $Y \setminus F$ may be embedded order-preservingly into $T$. As $\varphi \not\leqslant \mathrm{otp}(T)$, we have $\varphi \not\leqslant \mathrm{otp}(Y \setminus F)$. As $\varphi$ is unionwise indecomposable and infinite and $F$ is finite, it follows that $\varphi \not\leqslant \mathrm{otp}(Y) = \varphi$, a contradiction. □

The fact that $\mathfrak{b}$ has uncountable cofinality together with Corollary 3.2, Theorem 6.7 and Corollary 6.8 yields the following corollary summarising the central results of this paper:

**Corollary 6.18.** *Let $\varphi$ be a countable order type. If $\varphi$ is equimorphic to an order type in $\{0, 1, \omega^*, \omega, \eta\}$, then*

$$\binom{\mathfrak{b}}{\varphi} \longrightarrow \binom{\mathfrak{b} \quad \alpha}{\varphi \quad \varphi}$$

*for all $\alpha < \mathfrak{t}$; otherwise*

$$\binom{\mathfrak{b}}{\varphi} \not\longrightarrow \binom{\mathfrak{b} \quad 1}{\varphi \quad \varphi}.$$

## 7. Questions

We are interested in relations which can be proved using only ZFC. Therefore the following questions seem natural:

**Question B.** *Does the relation*

$$\begin{pmatrix} \varphi \\ \psi \end{pmatrix} \longrightarrow \begin{pmatrix} \varphi & n \\ \psi & n \end{pmatrix}$$

*hold for all countable unionwise indecomposable order types $\varphi$ and $\psi$ and all natural numbers $n$? For example, what about $\varphi = \omega^{\omega+1}$ while $\psi = \omega$ and $n = 2$?*

**Question C.** *Does the relation*

$$\begin{pmatrix} \omega_1 \\ \varphi \end{pmatrix} \longrightarrow \begin{pmatrix} \alpha & \alpha \\ \varphi & \varphi \end{pmatrix}$$

*necessarily hold for all countable ordinals $\alpha$ and all countable unionwise indecomposable order types $\varphi$?*

Note that if the scope of the question above is restricted by additionally requiring $\varphi$ to be an ordinal, the answer is positive, cf. [BH73, Corollary 3].

In light of Proposition 6.13, the following question suggests itself:

**Question D.** *Is it consistent that*

$$\begin{pmatrix} \omega_1 \\ \varphi \end{pmatrix} \longrightarrow \begin{pmatrix} \omega_1 & \alpha \\ \varphi & \varphi \end{pmatrix}$$

*for all countable ordinals $\alpha$ and all countable unionwise indecomposable order types $\varphi$?*

Furthermore, in light of Corollary 6.2 and Theorem 6.7, the following question arises:

**Question E.** *Does the relation*

$$\begin{pmatrix} \kappa \\ \omega \end{pmatrix} \longrightarrow \begin{pmatrix} \kappa & \alpha \\ \omega & \omega \end{pmatrix}$$

*hold for all uncountable cardinals $\kappa \leqslant \mathfrak{c}$ and all $\alpha < \min(\mathfrak{s}_{\aleph_0}, \operatorname{cf}(\kappa))$?*

Finally, the following question seems deceptively simple:

**Question F.** *Is $2$ the only countable linear order type which is unionwise decomposable yet multiplicatively untranscendable?*

Note that the answer is negative if "countable" is dropped from the question: In a conversation with the second author, Garrett Ervin pointed out that the order type $\lambda$ of the real numbers is unionwise decomposable yet multiplicatively untranscendable, as $\lambda$ is not embeddable into a product of two smaller types and $\lambda \nrightarrow (\lambda)^1_2$ (cf. [ER56, Theorem 9]).

## References

[Bel81]   Murray G. Bell, *On the Combinatorial Principle $P(\mathfrak{c})$*, Fund. Math. **114** (1981), no. 2, 149–157, DOI: 10.4064/FM-114-2-149-157, also available at https://eudml.org/doc/211293.

[BF11]    Jörg Brendle and Vera Fischer, *Mad Families, Splitting Families and Large Continuum*, J. Symb. Logic **76** (2011), no. 1, 198–208, DOI: 10.2178/JSL/1294170995, also available at http://www.logic.univie.ac.at/~vfischer/absplit.pdf.

[BH73]    James Earl Baumgartner and András Hajnal, *A Proof (Involving Martin's Axiom) of a Partition Relation*, Fund. Math. **78** (1973), no. 3, 193–203, DOI: 10.4064/FM-78-3-193-203.

[BJ95] Tomek Bartoszyński and Haim Judah, *Set Theory: On the Structure of the Real Line*, A K Peters, Wellesley, MA, 1995, DOI: 10.1112/S0024609396222374.

[Bla10] Andreas Blass, *Combinatorial Cardinal Characteristics of the Continuum*, Handbook of Set Theory (Matthew Foreman and Akihiro Kanamori, eds.), Springer, Dordrecht, 2010, pp. 395–489, DOI: 10.1007/978-1-4020-5764-9_7.

[BR14] Jörg Brendle and Dilip Raghavan, *Bounding, Splitting, and Almost Disjointness*, Ann. Pure Appl. Logic **165** (2014), no. 2, 631–651, DOI: 10.1016/J.APAL.2013.09.002, arXiv: 1211.5170 [math.LO].

[CGW20] William Chen, Shimon Garti, and Thilo Volker Weinert, *Cardinal Characteristics of the Continuum and Partitions*, Israel J. Math. **235** (2020), no. 1, 13–38, DOI: 10.1007/S11856-019-1942-y, arXiv: 1801.00238 [math.LO].

[CH17] Andrés Eduardo Caicedo and Jacob Hilton, *Topological Ramsey Numbers and Countable Ordinals*, Foundations of Mathematics, Contemp. Math., vol. 690, Amer. Math. Soc., Providence, RI, 2017, pp. 87–120, DOI: 10.1090/CONM/690/13864, arXiv: 1510.00078 [math.LO].

[EHM70] Paul Erdős, András Hajnal, and Eric Charles Milner, *Set Mappings and Polarized Partition Relations*, Combinatorial Theory and Its Applications, I (Proc. Colloq., Balatonfüred, 1969) (Paul Erdős, Alfréd Rényi, and Vera Turán Sós, eds.), North-Holland, Amsterdam, 1970, pp. 327–363, available at https://pdfs.semanticscholar.org/62df/2699160f4cc06e5968e66339e14a6d3911fa.pdf.

[EHM71] Paul Erdős, András Hajnal, and Eric Charles Milner, *Polarized Partition Relations for Ordinal Numbers*, Studies in Pure Mathematics (Leonid Mirsky, ed.), Academic Press, London, 1971, pp. 63–87, available at https://users.renyi.hu/~p_erdos/1971-17.pdf.

[EM72] Paul Erdős and Eric Charles Milner, *A Theorem in the Partition Calculus*, Canad. Math. Bull. **15** (1972), no. 4, 501–505, available at http://www.renyi.hu/~p_erdos/1972-03.pdf.

[ER56] Paul Erdős and Richard Rado, *A Partition Calculus in Set Theory*, Bull. Amer. Math. Soc. **62** (1956), 427–489, available at http://www.ams.org/journals/bull/1956-62-05/S0002-9904-1956-10036-0/S0002-9904-1956-10036-0.pdf.

[ES35] Paul Erdős and George Szekeres, *A Combinatorial Problem in Geometry*, Compos. Math. **2** (1935), 463–470, available at http://www.numdam.org/item/CM_1935__2__463_0/.

[FM17] Vera Fischer and Diego Alejandro Mejía, *Splitting, Bounding, and Almost Disjointness Can Be Quite Different*, Canad. J. Math. **69** (2017), no. 3, 502–531, DOI: 10.4153/CJM-2016-021-8, arXiv: 1508.01068 [math.LO].

[Fra00] Roland Fraïssé, *Theory of Relations*, Studies in Logic and the Foundations of Mathematics, vol. 145, North-Holland, Amsterdam, 2000.

[FS08] Vera Fischer and Juris Steprāns, *The Consistency of $\mathfrak{b} = \kappa$ and $\mathfrak{s} = \kappa^+$*, Fund. Math. **201** (2008), no. 3, 283–293, DOI: 10.4064/FM201-3-5.

[GL75] Frederick William Galvin and Jean Ann Larson, *Pinning Countable Ordinals*, Fund. Math. **82** (1975), no. 4, 357–361, DOI: 10.4064/FM-82-4-357-361.

[GS12] Shimon Garti and Saharon Shelah, *Combinatorial Aspects of the Splitting Number*, Ann. Comb. **16** (2012), no. 4, 709–717, DOI: 10.1007/S00026-012-0154-5, arXiv: 1007.2266 [math.LO].

[GS14] \_\_\_\_\_\_, *Partition Calculus and Cardinal Invariants*, J. Math. Soc. Japan **66** (2014), no. 2, 425–434, DOI: 10.2969/JMSJ/06620425, arXiv: 1112.5772 [math.LO].

[Hal17] Lorenz J. Halbeisen, *Combinatorial Set Theory: With a Gentle Introduction to Forcing*, Springer Monographs in Mathematics, Springer, London, $^2$2017, DOI: 10.1007/978-3-319-60231-8.

[Hau08] Felix Hausdorff, *Grundzüge einer Theorie der geordneten Mengen*, Math. Ann. **65** (1908), no. 4, 435–505, DOI: 10.1007/BF01451165.

[Hil16] Jacob Hilton, *The Topological Pigeonhole Principle for Ordinals*, J. Symb. Log. **81** (2016), no. 2, 662–686, DOI: 10.1017/JSL.2015.45, arXiv: 1410.2520 [math.LO].

[HS69] Labib Haddad and Gabriel Sabbagh, *Calcul de certains nombres de Ramsey généralisés*, C. R. Acad. Sci. Paris Sér. A-B **268** (1969), A1233–A1234.

[IRW] Ferdinand Ihringer, Deepak Rajendraprasad, and Thilo Volker Weinert, *New Bounds on the Ramsey Numbers $r(I_m, L_n)$*, arXiv: 1707.09556 [math.CO].

[Jon08] Albin Lee Jones, *On a Result of Szemerédi*, J. Symb. Logic **73** (2008), no. 3, 953–956, DOI: 10.2178/JSL/1230396758.

[KW96] Anastasis Kamburelis and Bogdan Zbigniew Węglorz, *Splittings*, Arch. Math. Logic **35** (1996), no. 4, 263–277, DOI: 10.1007/S001530050044.

[Lar73] Jean Ann Larson, *A Short Proof of a Partition Theorem for the Ordinal $\omega^\omega$*, Ann. Math. Logic **6** (1973), 129–145, DOI: 10.1016/0003-4843(73)90006-5.

[Lav71] Richard Laver, *On Fraïssé's Order Type Conjecture*, Ann. of Math. (2) **93** (1971), 89–111, available at https://www.jstor.org/stable/pdf/1970754.pdf.

[LHW19] Chris Lambie-Hanson and Thilo Volker Weinert, *Partitioning Subsets of Generalised Scattered Orders*, J. Math. Soc. Japan **71** (2019), no. 1, 235–257, DOI: 10.2969/JMSJ/78617861, arXiv: 1701.05791 [math.LO].

[Mal89] Vyacheslav Ivanovich Malyhin, *Topological Properties of Cohen Generic Extensions*, Trans. Mosc. Math. Soc. **52** (1989), 3–33, available at http://www.mathnet.ru/links/cc672532bfb79cbe49268c92dbe18bbe/mmo483.pdf.

[Mer19] Omer Mermelstein, *Calculating the Closed Ordinal Ramsey Number $R^{cl}(\omega \cdot 2, 3)$*, Israel J. Math. **230** (2019), no. 1, 387–407, DOI: 10.1007/S11856-019-1827-0, arXiv: 1702.03878 [math.LO].

[MS13] Maryanthe Malliaris and Saharon Shelah, *General Topology Meets Model Theory, on $\mathfrak{p}$ and $\mathfrak{t}$*, Proc. Natl. Acad. Sci. USA **110** (2013), no. 33, 13300–13305, DOI: 10.1073/PNAS.1306114110.

[Orr95] John Lindsay Orr, *Shuffling of Linear Orders*, Canad. Math. Bull. **38** (1995), no. 2, 223–229, DOI: 10.4153/CMB-1995-032-1.

[Ros82] Joseph Goeffrey Rosenstein, *Linear Orderings*, Pure and Applied Mathematics, vol. 98, Academic Press, New York, 1982.

[RT18] Dilip Raghavan and Stevo Todorcevic, *Suslin Trees, the Bounding Number, and Partition Relations*, Israel J. Math. **225** (2018), no. 2, 771–796, DOI: 10.1007/S11856-018-1677-1, arXiv: 1602.07901 [math.LO].

[Sch10] Rene Schipperus, *Countable Partition Ordinals*, Ann. Pure Appl. Logic **161** (2010), no. 10, 1195–1215, DOI: 10.1016/J.APAL.2009.12.007.

[Sch12] \_\_\_\_\_\_, *The Topological Baumgartner–Hajnal Theorem*, Trans. Amer. Math. Soc. **364** (2012), no. 8, 3903–3914, DOI: 10.1090/S0002-9947-2012-04990-7.

[Ste93] Juris Steprāns, *Combinatorial Consequences of Adding Cohen Reals*, Set Theory of the Reals (Haim Judah, ed.), Israel Mathematics Conference Proceedings, vol. 6, Gelbart Research Institute for Mathematical Sciences, Bar-Ilan University, Ramat-Gan, 1993, pp. 583–617.

[vD84] Eric Karel van Douwen, *The Integers and Topology*, Handbook of Set-Theoretic Topology (Kenneth Kunen and Jerry E. Vaughan, eds.), North-Holland, Amsterdam, 1984, pp. 111–167, DOI: 10.1016/B978-0-444-86580-9.50006-9.

[Wei14] Thilo Volker Weinert, *Idiosynchromatic Poetry*, Combinatorica **34** (2014), no. 6, 707–742, DOI: 10.1007/S00493-011-2980-1, also available at https://www.math.bgu.ac.il/~weinert/Poetry_colour0.pdf.

[Wil77] Neil Hale Williams, *Combinatorial Set Theory*, Studies in Logic and the Foundations of Mathematics, vol. 91, North-Holland, Amsterdam, 1977.

Institute of Discrete Mathematics and Geometry, TU Wien, Wiedner Hauptstrasse 8–10/104, 1040 Wien, Austria
*Email address*: `mail@l17r.eu`
*URL*: `https://l17r.eu`

Kurt Gödel Research Center for Mathematical Logic, University of Vienna, Währinger Strasse 25, 1090 Wien, Austria
*Email address*: `thilo.weinert@univie.ac.at`